\documentclass[11pt]{article}
\usepackage{amsfonts}
\usepackage{bbm}
\usepackage{mathrsfs}
\usepackage{amssymb}
\usepackage{url}

\title{ On The Algebras $U_q^{\pm}(A_N)$:\\ From A Constructive-Computational Viewpoint}

\author{Rabigul Tuniyaz\thanks{Project supported by
the National Natural Science Foundation of China (1186106).}\\
{\small Department of Mathematics, School of  Science}\\
{\small Xinjiang Institute of Science and Technology}\\
{\small  Akesu, 843100, Xinjiang, China}\\}

\date{}

\begin{document}
\maketitle
\begin{center}
\begin{minipage}{120mm}
{\small {\bf Abstract.} Let $U_q^+(A_N)$ (resp.  $U_q^-(A_N)$) be the $(+)$-part (resp. $(-)$-part) of the Drinfeld-Jimbo quantum group of type $A_N$ over a field $K$. With respect to  Jimbo relations and the PBW $K$-basis ${\cal B}$ of $U_q^+(A_N)$ (resp.  $U_q^-(A_N)$) established by Yamane, it is shown, by constructing an appropriate monomial ordering $\prec$ on ${\cal B}$,  that $U_q^+(A_N)$ (resp.  $U_q^-(A_N)$) is a solvable polynomial algebra. Consequently,  further structural properties of $U_q^+(A_N)$ (resp.  $U_q^-(A_N)$) and their modules may be established and realized in a constructive-computational way.}
\end{minipage}\end{center} {\parindent=0pt\par

{\bf Key words:} Quantum group, Solvable polynomial algebra, Gr\"obner basis}
\vskip -.5truecm

\renewcommand{\thefootnote}{\fnsymbol{footnote}}
\let\footnote\relax\footnotetext{E-mail: rabigul802@sina.com}
\let\footnote\relax\footnotetext{2010 Mathematics Subject Classification: 16T20, 16Z05.}

\def\NZ{\mathbb{N}}
\def\QED{\hfill{$\Box$}}
\def \r{\rightarrow}
\def\mapright#1#2{\smash{\mathop{\longrightarrow}\limits^{#1}_{#2}}}

\def\v5{\vskip .5truecm}
\def\OV#1{\overline {#1}}
\def\hang{\hangindent\parindent}
\def\textindent#1{\indent\llap{#1\enspace}\ignorespaces}
\def\item{\par\hang\textindent}

\def\LH{{\bf LH}}\def\LM{{\bf LM}}\def\LT{{\bf
LT}}\def\KX{K\langle X\rangle} \def\KS{K\langle X\rangle}
\def\B{{\cal B}} \def\LC{{\bf LC}} \def\G{{\cal G}} \def\FRAC#1#2{\displaystyle{\frac{#1}{#2}}}
\def\SUM^#1_#2{\displaystyle{\sum^{#1}_{#2}}}
\def\T{\widetilde}\def\KD{K\langle D\rangle}
\def\PRC{\prec_{d\textrm{\tiny -}lex}}

\section*{1. Introduction}

Let $K$ be a field, and let $U_q(A_N)$ be the Drinfeld-Jimbo quantum group of type $A_N$ in the sense of ([Dr], [Jim]), that is, with $q\in K$, $q^8\ne 1$, and a positive integer $N$, $U_q(A_N)$ is generated by $\{E_{i},{K_{i}}^{\pm1},F_{i}|1\leq i,j\leq N\}$ over $K$  subject to the relations
$$\begin{array}{rcl}
K&=&\left\{\begin{array}{l} K_{i}K_{j}-K_{j}K_{i},~ K_{i}{K_{i}}^{-1}-1, ~{K_{i}}^{-1}K_{i}-1,\\
E_{j}{K_{i}}^{\pm 1}-q^{{\mp d_{i}}a_{ij}}{K_{i}}^{\pm 1}E_{j},~
K_{i}^{\pm 1}F_{j}-q^{{\mp d_{i}}a_{ij}}{F_{j}K_{i}}^{\pm 1}\end{array}\right\}\\ \\
T&=&\left\{E_{i}F_{j}-F_{j}E_{i}-\delta_{ij}\frac{K_{i}^{2}-K_{i}^{-2}}{q^{2d_{i}}
-q^{-2d_{i}}}\right\};\\
\\
S^{+}&=&\left\{\left.\displaystyle{\sum_{v=0}^{1-a_{ij}}}(-1)^v
\left[ \begin{array}{c}1-a_{ij} \\ v\end{array} \right]_tE_{i}^{1-a_{ij}-v}E_{j}E_{i}^{v}~\right |~i\ne j, t=q^{2d_{i}}\right\};
\end{array}$$
$$\begin{array}{rcl}
S^{-}&=&\left\{\left.\displaystyle{\sum_{v=0}^{1-a_{ij}}}(-1)^v
\left[ \begin{array}{c}1-a_{ij} \\v\end{array} \right]_{t}F_{i}^{1-a_{ij}-v}F_{j}F_{i}^v~\right |~i\ne j, t=q^{2d_{i}}\right\},
\end{array}$$
where
$$\left[ \begin{array}{c}
m\\
n\end{array} \right]_{\alpha}=\left\{\begin{array}{ll}
\prod_{i=1}^{n}\frac{t^{m-1+i}-t^{i-m-1}}{t^{i}-t^{-i}},&m>n>0,\\
1&n=0~\hbox{or}~n=m.\end{array}\right.$$
Let $U_{q}^0(A_N)$, $U_{q}^+(A_N)$, and
$U_{q}^-(A_N)$ be the subalgebras of $U_q(A_N)$ generated by
$\{{K_i}^{\pm1}~|~1\leq i\leq N\}$, $\{E_i~|~ 1\leq i\leq N\}$, and
$\{F_i~|~ 1\leq i\leq N\}$, respectively. Then $U_{q}(A_N)$ has the
triangular decomposition
$$U_q(A_N)\cong{U_{q}^+(A_N)}\otimes_K {U_q^0(A_N)}\otimes_K {U_q^-(A_N)},$$
where $U_q^+(A_N)$ and $U_q^-(A_N)$ are  called the $(+)$-part and ($-$)-part of $U_q(A_N)$ respectively. In [Ros] and [yam], it was proved that with respect to the Jimbo defining relations, $U^+_q(A_N)$ (similarly $U_q^-(A_N)$) has the standard  PBW $K$-basis $\B$ (see Section 2 for this basis). In this paper, we show (in Section 2) that there is a monomial ordering $\prec$ on the PBW basis $\B$ of  $U_q^+(A_N)$ (resp.  $U_q^-(A_N)$) that makes $U_q^+(A_N)$ (resp.  $U_q^-(A_N)$) into a solvable polynomial algebra in the sense of [K-RW]. In Section 3, we show that the main result (Theorem 2.3) obtained in Section 2 may  enable us  to  establish and realize further structural properties of $U_q^+(A_N)$ (resp.  $U_q^-(A_N)$) and their modules in a constructive-computational way.\v5

Throughout this note, $K$ denotes a field of characteristic 0, $K^*=K-\{0\}$, and all $K$-algebras considered are associative with multiplicative identity 1. If $S$ is a nonempty subset of an algebra $A$, then we write $\langle S\rangle$ for the two-sided ideal of $A$ generated by $S$.\v5

\section*{2. $U_q^+(A_N)$ (resp.  $U_q^-(A_N)$) is a solvable polynomial algebra}\par
We start  by recalling from ([K-RW], [Li1, 6]) the following definitions and notations. Suppose that a finitely generated  $K$-algebra $A=K[a_1,\ldots ,a_n]$ has the PBW $K$-basis $\B =\{ a^{\alpha}=a_{1}^{\alpha_1}\cdots
a_{n}^{\alpha_n}~|~\alpha =(\alpha_1,\ldots ,\alpha_n)\in\NZ^n\}$, and that $\prec$ is a total ordering on $\B$. Then every nonzero element $f\in A$ has a unique expression
$$\begin{array}{rcl} f&=&\lambda_1a^{\alpha (1)}+\lambda_2a^{\alpha (2)}+\cdots +\lambda_ma^{\alpha (m)},\\
&{~}&\hbox{such that}~a^{\alpha (1)}\prec a^{\alpha
(2)}\prec\cdots \prec a^{\alpha (m)},\\
&{~}&\hbox{where}~ \lambda_j\in K^*,~a^{\alpha
(j)}=a_1^{\alpha_{1j}}a_2^{\alpha_{2j}}\cdots a_n^{\alpha_{nj}}\in\B
,~1\le j\le m.
\end{array}$$
Since elements of $\B$ are conventionally called {\it monomials}\index{monomial}, the {\it leading monomial of $f$} is defined as $\LM
(f)=a^{\alpha (m)}$, the {\it leading coefficient of $f$} is defined
as $\LC (f)=\lambda_m$, and the {\it leading term of $f$} is defined
as $\LT (f)=\lambda_ma^{\alpha (m)}$.\v5

{\bf Definition 2.1}  Suppose that the $K$-algebra
$A=K[a_1,\ldots ,a_n]$ has the PBW basis $\B$. If $\prec$ is a
total ordering on $\B$ that satisfies the following three
conditions:{\parindent=1.35truecm\par

\item{(1)} $\prec$ is a well-ordering (i.e., every nonempty subset of $\B$ has a minimal element);\par

\item{(2)} For $a^{\gamma},a^{\alpha},a^{\beta},a^{\eta}\in\B$, if $a^{\gamma}\ne 1$, $a^{\beta}\ne
a^{\gamma}$, and $a^{\gamma}=\LM (a^{\alpha}a^{\beta}a^{\eta})$,
then $a^{\beta}\prec a^{\gamma}$ (thereby $1\prec a^{\gamma}$ for
all $a^{\gamma}\ne 1$);\par

\item{(3)} For $a^{\gamma},a^{\alpha},a^{\beta}, a^{\eta}\in\B$, if
$a^{\alpha}\prec a^{\beta}$, $\LM (a^{\gamma}a^{\alpha}a^{\eta})\ne
0$, and $\LM (a^{\gamma}a^{\beta}a^{\eta})\not\in \{ 0,1\}$, then
$\LM (a^{\gamma}a^{\alpha}a^{\eta})\prec\LM
(a^{\gamma}a^{\beta}a^{\eta})$,\par}{\parindent=0pt
then $\prec$ is called a {\it monomial ordering} on $\B$ (or a
monomial ordering on $A$).} \v5

{\bf Definition 2.2} A finitely generated $K$-algebra $A=K[a_1,\ldots ,a_n]$
is called a {\it solvable polynomial algebra} if $A$ has the PBW $K$-basis $\B =\{
a^{\alpha}=a_1^{\alpha_1}\cdots a_n^{\alpha_n}~|~\alpha
=(\alpha_1,\ldots ,\alpha_n)\in\NZ^n\}$ and a monomial ordering $\prec$ on $\B$, such that
for $\lambda_{ji}\in K^*$ and  $f_{ji}\in A$,
$$\begin{array}{l} a_ja_i=\lambda_{ji}a_ia_j+f_{ji},~1\le i<j\le n,\\
\LM (f_{ji})\prec a_ia_j~\hbox{whenever}~f_{ji}\ne 0.\end{array}$$\par

It follows from [K-RW] that a solvable polynomial algebra $A$ is equipped with an algorithmic Gr\"obner basis theory, that is, every (two-sided, respectively one-sided) ideal of $A$ and every submodule of a free (left) $A$-module has a finite Gr\"obner basis which can be produced by running a noncommutative Buchberger Algorithm with respect to a given monomial ordering. It is also well known that  nowadays the noncommutative Buchberger Algorithm  for solvable polynomial algebras and their modules has been successfully implemented in the computer algebra system \textsf{Plural} [LS].  Concerning  basic constructive-computational theory and methods for solvable polynomial algebras and their modules, one is referred to [Li6] for more details.\v5

Now, we aim to prove the following result.\v5

{\bf Theorem 2.3} If $q^8\ne 1$, then the algebra $U_q^+(A_N)$ is a solvable polynomial algebra in the sense of Definition 2.2. \vskip 6pt

{\bf Proof} First, recall that the Jimbo relations, namely the defining relations of $U_q^+(A_N)$ as described in [Yam], are given by
$$\begin{array}{ll} f_{13}=x_{mn}x_{ij}-q^{-2}x_{ij}x_{mn},&((i,j),(m,n))\in C_1\cup C_3,\\
f_{26}=x_{mn}x_{ij}-x_{ij}x_{mn},&((i,j),(m,n))\in C_2\cup C_6,\\
f_4=x_{mn}x_{ij}-x_{ij}x_{mn}+(q^2-q^{-2})x_{in}x_{mj},&((i,j),(m,n))\in
C_4,\\
f_5=x_{mn}x_{ij}-q^2x_{ij}x_{mn}+qx_{in},&((i,j),(m,n))\in C_5,\end{array}$$
where with $\Lambda_N=\{ (i,j)\in\NZ\times
\NZ~|~1\le i<j\le N+1\}$, the $C_i$ are given by
$$\begin{array}{l} C_1=\{ ((i,j),(m,n))~|~i=m<j<n\},\\
C_2=\{ ((i,j),(m,n))~|~i<m<n<j\},\\
C_3=\{ ((i,j),(m,n))~|~i<m<j=n\},\\
C_4=\{ ((i,j),(m,n))~|~i<m<j<n\},\\
C_5=\{ ((i,j),(m,n))~|~i<j=m<n\},\\
C_6=\{((i,j),(m,n))~|~i<j<m<n\}.\end{array}$$
It follows from [Yam] that for $q^8\ne
1$, $U^+_q(A_N)$ has the standard PBW $K$-basis
$$\B =\left\{1,~x_{i_1j_1}x_{i_2j_2}\cdots
x_{i_kj_k}~\left | ~(i_{\ell},j_{\ell})\in\Lambda_N,~k\ge 1,~
(i_1,j_1)\le_{lex} (i_2,j_2)\le_{lex}\cdots \le_{lex} (i_k,j_k)\right.\right\} ,$$
where $<_{lex}$ is the lexicographic ordering on $\Lambda_N$, i.e.,
$$(l,k)<_{lex}(i,j)\Leftrightarrow\left\{\begin{array}{l} l<i,\\
\hbox{or}~l=i~\hbox{and}~k<j.\end{array}\right.$$\par

We now start on constructing a monomial ordering on $\B$. In doing so, we let $X=\{ x_{ij}~|~(i,j)\in\Lambda_N\}$ and introduce an ordering $\prec$ on $X$: for $x_{lk}$, $x_{ij}\in X$,
$$x_{lk}\prec x_{ij}\Leftrightarrow\left\{\begin{array}{l}
l<i,\\
\hbox{or}~l=i~\hbox{and}~k>j.\end{array}\right.$$
Note that the ordering $\prec$ is not the one introduced by the lexicographic ordering $<_{lex}$ on $\Lambda_N$. Furthermore, we extend $\prec$ to $\B$:
$$1\prec u~\hbox{for all}~u=x_{i_1j_1}x_{i_2j_2}\cdots x_{i_rj_r}\in\B -\{ 1\},$$
and for $u=x_{i_1j_1}x_{i_2j_2}\cdots x_{i_rj_r}$,
$v=x_{l_1t_1}x_{l_2t_2}\cdots x_{l_ht_h}\in\B$,
$$u\prec v\Leftrightarrow\left\{\begin{array}{l} r<h~\hbox{and}~
x_{i_1j_1}=x_{l_1t_1},~x_{i_2j_2}=x_{l_2t_2},\ldots ,x_{i_rj_r}=x_{l_rt_r},\\
\hbox{or there exists an}~m, 1\le m\le r,~\hbox{such that}\\
x_{i_1j_1}=x_{l_1t_1},~x_{i_2j_2}=x_{l_2t_2},\ldots ,x_{i_{m-1}j_{m-1}}=x_{l_{m-1}t_{m-1}}\\
\hbox{but}~x_{i_mj_m}\prec x_{l_mt_m}.\end{array}\right.$$
It is straightforward to check that $\prec$ is reflexive, antisymmetrical, transitive, and any two generators $x_{ij}, x_{kl}\in X$ are comparable, thereby $\prec$ is a
total ordering on $\B$. Also since $\Lambda_N$ is a finite set, it can be directly verified  that $\prec$ satisfies the descending chain condition on $\B$, namely $\prec$ is a well-ordering on $\B$. \par
It remains to show that $\prec$ satisfies the conditions (2) and (3) of Definition 2.1, and that with respect to $\prec$ on $\B$, the relations $f_{13}$, $f_{26}$, $f_4$ and $f_5$satisfied by generators of $U_q^+(A_N)$ (which are given by the Jimbo relations) have the property required by Definition 2.2. To see this, let $x_{mn}, x_{kl}, x_{ij}\in X$, and suppose that $x_{mn}\prec x_{kl}$. If $((i,j), (m,n))\in C_4$, then since
$i<m<j<n$ and $x_{in}x_{mj}\in \B$,  the jimbo relation $f_4$ gives rise to
$$\begin{array}{rcl} x_{mn}x_{ij}&=&x_{ij}x_{mn}-(q^2-a^{-2})x_{in}x_{mj}\\
&{~}&\hbox{with}
~\LM ((a^2-a^{-2})x_{in}x_{mj})=x_{in}x_{mj}\prec x_{ij}x_{mn}=\LM (x_{mn}x_{ij}).
\end{array}$$
On the other hand, since $i<j$, if $((i,j), (k,l))\in C_4$, then noticing  $i<k<j<l$ and $x_{il}x_{kj}\in\B$, the jimbo relation $f_4$ gives rise to
$$\begin{array}{rcl} x_{kl}x_{ij}&=&x_{ij}x_{kl}-(q^2-a^{-2})x_{il}x_{kj}\\
&{~}&\hbox{with}
~\LM ((a^2-a^{-2})x_{il}x_{kj})=x_{il}x_{kj}\prec x_{ij}x_{kl}=\LM (x_{kl}x_{ij}).
\end{array}$$
Thus, we have shown that  if $((i,j), (m,n)), ((i,j), (k,l))\in C_4$, then
$$\begin{array}{l} x_{mn}\prec x_{kl}~\hbox{implies}~\LM (x_{mn}x_{ij})=x_{ij}x_{mn}\prec x_{ij}x_{kl}=\LM (x_{kl}x_{ij}),\\
\hbox{and the generating relations of}~U_q^+(A_N)~\hbox{determined by}~f_4\\
\hbox{have the property reqired by Definition 2.2}.\end{array}\eqno{(1)}$$
Similarly in the case that $((m,n), (i,j)), ((k,l), (i,j))\in C_4$, we  have
$$\begin{array}{l} x_{mn}\prec x_{kl}~\hbox{implies}~\LM (x_{ij}x_{mn})=x_{mn}x_{ij}\prec x_{kl}x_{ij}=\LM (x_{ij}x_{kl}),\\
\hbox{and the generating relations of}~U_q^+(A_N)~\hbox{determined by}~f_4\\
\hbox{have the property reqired by Definition 2.2}.\end{array}\eqno{(2)}$$
Fuethermore, if  $((i,j), (m,n)), ((i,j), (k,l))\in C_5$, then since $i<j=m<n$ and $x_{in}, x_{il}\in\B$, the jimbo relation $f_5$ gives rise to
$$\begin{array}{rcl} x_{mn}x_{ij}&=&q^2x_{ij}x_{mn}-qx_{in}\\
&{~}&\hbox{with}
~\LM (qx_{in})=x_{in}\prec x_{ij}x_{mn}=\LM (x_{mn}x_{ij}),\\
x_{kl}x_{ij}&=&q^2x_{ij}x_{kl}-qx_{il}\\
&{~}&\hbox{with}
~\LM (qx_{il})=x_{il}\prec x_{ij}x_{kl}=\LM (x_{kl}x_{ij}).
\end{array}$$
This shows that  if $((i,j), (m,n)), ((i,j), (k,l))\in C_5$, then
$$\begin{array}{l} x_{mn}\prec x_{kl}~\hbox{implies}~\LM (x_{mn}x_{ij})=x_{ij}x_{mn}\prec x_{ij}x_{kl}=\LM (x_{kl}x_{ij}),\\
\hbox{and the generating relations of}~U_q^+(A_N)~\hbox{determined by}~f_5\\
\hbox{have the property reqired by Definition 2.2},\end{array}\eqno{(3)}$$
and in the case that $((m,n), (i,j)), ((k,l), (i,j))\in C_5$, we also have
$$\begin{array}{l} x_{mn}\prec x_{kl}~\hbox{implies}~\LM (x_{ij}x_{mn})=x_{mn}x_{ij}\prec x_{kl}x_{ij}=\LM (x_{ij}x_{kl}),\\
\hbox{and the generating relations of}~U_q^+(A_N)~\hbox{determined by}~f_5\\
\hbox{have the property reqired by Definition 2.2}.\end{array}\eqno{(4)}$$
At this stage, the relations $f_{13}$, $f_{26}$, $f_4$, and $f_5$, the  assertions $(1)$, $(2)$, $(3)$, and $(4)$ derived above, all together enable us to conclude that for any $x_{mn}$, $x_{kl}$, $x_{ij}\in X$, if $x_{mn}\prec x_{kl}$ then
then
$$\begin{array}{l} x_{mn}\prec x_{kl}~\hbox{implies}~\LM (x_{ij}x_{mn})=x_{mn}x_{ij}\prec x_{kl}x_{ij}\LM (x_{ij}x_{kl}),\\
x_{mn}\prec x_{kl}~\hbox{implies}~\LM (x_{mn}x_{ij})=x_{ij}x_{mn}\prec x_{ij}x_{kl}=\LM (x_{kl}x_{ij}),\\
\hbox{and the generating relations of}~U_q^+(A_N)~\hbox{determined by}~f_{13},f_{26}, f_4,~\hbox{and}~ f_5\\
\hbox{have the property reqired by Definition 2.2}.\\
\end{array}\eqno{(5)}$$
Finally, by means of the assertions $(1)$, $(2)$, $(3)$, $(4)$, and $(5)$ derived above, it is straightforward to check that the conditions (2) and (3) of Definition 2.1 are satisfied by $\prec$, thereby $\prec$ is a monomial ordering on $\B$, and consequently, $U_q^+(A_N)$ is a solvable polynomial algebra n the sense of Definition 2.2,  as desired.\QED\v5

Similarly, the following assertion holds.\v5

{\bf Theorem 2.4} Let $U_q^-(A_N)$ be the $(-)$-part of the Drinfeld-Jimbo quantum group of type $A_N$. Then $U_q^-(A_N)$ is a solvable polynomial algebra in the sense of Definition 2.2.\par\QED\v5

Also by [L6, Proposition 1.1.6] we have the following\v5

{\bf Theorem 2.5} The tensor product $R=U_q^+(A_N)\otimes_KU_q^-(A_N)$ is a solvable polynomial algebra, where, for convenience, if $\B_1$ and $\B_2$ denote the PBW bases of  $U_q^+(A_N)$ and $U_q^-(A_N)$ respectively, $\prec_1$ and $\prec_2$ denote the monomial orderings of $U_q^+(A_N)$ and $U_q^-(A_N)$ respectively,
then $\B=\{u\otimes
v~|~u\in\B_1,~v\in\B_2\}$ is a PBW basis of $R$, and a monomial ordering $\prec$ on $\B$ is
defined subject to the rule: for $u_1\otimes v_1$, $u_2\otimes v_2\in\B$,
$$u_1\otimes
v_1\prec u_2\otimes
v_2\Leftrightarrow\left\{\begin{array}{l}
u_1\prec_1 u_2\\
\hbox{or}\\
u_1=u_2~\hbox{and}~v_1\prec_2 v_2.\end{array}\right.$$

\section*{3. Applications of Theorem 2.3 and Theorem 2.4}\par

Let $U_q^+(A_N)$ (resp. $U^-_q(A_N$) be the $(+)$-part (resp. $(-)$-part) of the Drinfeld-Jimbo quantum group of type $A_N$. In this section, we show that Theorem 2.3 and Theorem 2.4 obtained in the last section may enable us  to establish and realize further   structural properties of the algebra $U_q^+(A_N)$ (resp. $U^-_q(A_N)$ and their modules in a constructive-computational way. For more details on the basic constructive-computational theory and methods for solvable polynomial algebras and their modules, one is referred to [Li6]. All notions and notations used in previous sections are retained.\v5

As the first application of Theorem 2.3, we  recapture a known result (see [Li3, P.135, Example 2]) concerning the defining relations of $U_q^+(A_N)$.\v5

{\bf Proposition 3.1} Let $\KX =K\langle X_{ij}~|~(i,j)\in\Lambda_N\rangle$ be the free $K$-algebra generated by $X=\{X_{ij}~|~(i,j)\in\Lambda_N\}$, and
$$\G =\left\{\begin{array}{ll} F_{13}=X_{mn}X_{ij}-q^{-2}X_{ij}X_{mn},&((i,j),(m,n))\in C_1\cup C_3,\\
F_{26}=X_{mn}X_{ij}-X_{ij}X_{mn},&((i,j),(m,n))\in C_2\cup C_6,\\
F_4=X_{mn}X_{ij}-X_{ij}X_{mn}+(q^2-q^{-2})X_{in}X_{mj},&((i,j),(m,n))\in
C_4,\\
F_5=X_{mn}X_{ij}-q^2X_{ij}X_{mn}+qX_{in},&((i,j),(m,n))\in C_5,\end{array}\right\}$$
the set of defining relations of $U_q^+(A_N)$, where with $\Lambda_N=\{ (i,j)\in\NZ\times
\NZ~|~1\le i<j\le N+1\}$, the $C_i$s are given by
$$\begin{array}{l} C_1=\{ ((i,j),(m,n))~|~i=m<j<n\},\\
C_2=\{ ((i,j),(m,n))~|~i<m<n<j\},\\
C_3=\{ ((i,j),(m,n))~|~i<m<j=n\},\\
C_4=\{ ((i,j),(m,n))~|~i<m<j<n\},\\
C_5=\{ ((i,j),(m,n))~|~i<j=m<n\},\\
C_6=\{((i,j),(m,n))~|~i<j<m<n\}.\end{array}$$
Then, there exists a monomial ordering $\prec_{_X}$ on $\KX$ such that $\G$ is a Gr\"obner basis (in the sense of [Mor]) of the ideal $\langle\G\rangle$, and
$$\LM (\G)=\{ X_{mn}X_{ij}|~((i,j), (m,n))\in C_1\cup C_2\cup\cdots\cup C_6,~(i,j)<_{lex} (m,n)\}.$$
where $\LM (\G )$ is the set of leading monomials $\LM (g)$ of elements $g\in\G$ (for an $F\in\KX$, $\LM (F)$ is defined with respect to $\prec_{_X}$, as that defined for an element in a solvable polynomial algebra in the last section).\par
A similar result holds true for $U_q^-(A_N)$.\vskip 6pt

{\bf Proof} Since $U_q^+(A_N)$ is a solvable polynomial algebra by Theorem 2.3, this follows from a constructive characterization of solvable polynomial algebras [Li4, Theorem 2.1] (see also [Li6, Theorem 1.2.1]).\QED\v5

The next proposition stems from [Li2, Section 6, Example 1] and [Li3, P.167 Example 3; Ch.7, Section 5.7, Corollary 7.6; section 6.3, Corollary 3.2]. \v5

{\bf Proposition 3.2} [Li6, Proposition A1.16, Proposition 1.2.2]  Let $A=K[a_1,\ldots, a_n]$ be a solvable polynomial algebra with admissible system $(\B ,\prec )$. Then $A\cong \KX /\langle\G\rangle$, where $\KX =K\langle X_1,\ldots ,X_n\rangle$ is the free $K$-algebra of $n$ generators and $\G$ is a Gr\"obner basis of the ideal $\langle\G\rangle$ in $\KX$ with respect to some monomial ordering $\prec$ such that $\LM (\G )=\{ X_jX_i~|~1\le i<j\le n\}$, and the following statements hold.\par

(i) $A$ has Gelfand-Kirillov dimension GK.dim$A=n$.\par

(ii)  $A$ has global (homological) dimension gl.dim$A\le n$. If $\KX$ is  equipped
with a positive-weight $\NZ$-gradation and $\G$ is a homogeneous Gr\"obner basis, then gl.dim$A=n$.\par

(iii)   If $\KX$ is $\NZ$-graded by assigning each $X_i$ the degree 1, and $\G$ is a
homogeneous Gr\"obner basis, then $A$ is a homogeneous $2$-Koszul algebra;
otherwise, $A$ is a non-homogeneous $2$-Koszul algebra in the sense of ([Pos], [BG]).\par\QED

For the sake of saving notation, except for retaining all notions and other notations used in previous sections, in what follows we use  $R^+$ (resp. $R^-$) to denote  the algebra $U_q^+(A_N)$ (resp. $U_q^-(A_N)$).\v5

{\bf Theorem 3.3} With notation fixed above, the following statements hold.\par
(i) $R^+$ (resp. $R^-$) is a Noetherian domain.\par
(ii) $R^+$ (resp. $R^-$) has Gelfand-Kirillov dimension $\frac{N(N+1)}{2}$. \par
(iii) $R^+$ (resp. $R^-$) has global homological dimension $\le \frac{N(N+1)}{2}$.\par
(iv) $R^+$ (resp. $R^-$) is a non-homogeneous 2-Koszul algebra in the sense of ([Pos], [BG2]).\vskip 6pt

{\bf Proof} We prove all assertions for $R^+$, because similar argumentation works well for $R^-$.\par

(i) Though this result have been known from the literature (see[Yam]), here we emphasize that this property may follow immediately from Theorem 2.3. More precisely, that $R^+$ has no divisors of zero follows from the fact that $\LM (fg)=\LM (f)\LM (g)$ for all nonzero $f,g\in R^+$, and that the Noetherianess of $R^+$ follows from the fact that every nonzero one-sided ideal has a finite Gr\"obner basis (see [K-RW]). \par

Note that $R^+$ is a solvable polynomial algebra by Theorem 2.3, and it has $\frac{N(N+1)}{2}$ generators $x_{ij}$ (see the proof of Theorem 2.3). The assertions (ii), (iii), and (iv) follow from Proposition 3.1 and Proposition 3.2 above. \QED\v5

Before stating the next result,  let us recall the Auslander regularity and the Cohen-Macaulay property of an algebra for the reader's convenience.
A finitely generated algebra $A$ is said to \par
(a) be {\it Auslander regular} if $A$ has finite global homological dimension, and for every finitely generated left $A$-module $M$, every integer $j\ge 0$, and every (right) $A$-submodule $N$ of Ext$^j_A(M,A)$ we have that $j(N)\ge j$, where $j(N)$ is the grade number of $N$ which is the least integer $i$ such that Ext$^i_A(M,A)\ne 0$;\par

(b) satisfy the {\it Cohen-Macaulay property} if for every finitely generated left $A$-module $M$ we have the equality: GK.dim$M+j(M)=$ GK.dim$A$, where GK.dim denotes the Gelfand-Kirollov dimension of a module. \par

Concerning the Auslander regularity and the Cohen-Macaulay property of an algebra, particularly, an algebra  with filtered-graded structures, one is referred to [LVO].\v5

{\bf Theorem 3.4}  With notation as above, the following statements hold.\par
(i) $R^+$ (resp. $R^-$) is an Auslander regular algebra satisfying the Cohen-Macaulay property.\par
(ii) The $K_0$-group of $R^+$ (resp. $R^-$) is isomorphic to $\mathbb{Z}$, the additive group of integers.\vskip 6pt

{\bf Proof} We prove the two assertions only for $R^+$, because similar argumentation works well for $R^-$.\par
(i) Our approach is to employ certain specified filtered-graded structures associated with $R^+$.  As in the proof of Theorem 2.3, let $\Lambda_N=\{ (i,j)\in\NZ\times
\NZ~|~1\le i<j\le N+1\}$, $X=\{ x_{ij}~|~(i,j)\in\Lambda_N\}$,
$$\begin{array}{l} C_1=\{ ((i,j),(m,n))~|~i=m<j<n\},\\
C_2=\{ ((i,j),(m,n))~|~i<m<n<j\},\\
C_3=\{ ((i,j),(m,n))~|~i<m<j=n\},\\
C_4=\{ ((i,j),(m,n))~|~i<m<j<n\},\\
C_5=\{ ((i,j),(m,n))~|~i<j=m<n\},\\
C_6=\{((i,j),(m,n))~|~i<j<m<n\},\end{array}$$
and $$\B =\left\{ 1,~x_{i_1j_1}x_{i_2j_2}\cdots
x_{i_kj_k}~\left | ~(i_{\ell},j_{\ell})\in\Lambda_N,~k\ge 1,~
(i_1,j_1)\le_{lex} (i_2,j_2)\le_{lex}\cdots \le_{lex} (i_k,j_k)\right.\right\}$$
which is the PBW $K$-basis of $R^+$. Furthermore,  for every $x_{ij}\in X$, we assign the degree $d(x_{ij})=1$, so that each standard monomial $u=x_{i_1j_1}x_{i_2j_2}\cdots x_{i_kj_k}\in\B$ has a unique degree $$d(u)=d(x_{i_1j_1})+d(x_{i_2j_2})+\cdots +d(x_{i_kj_k}).$$
Now, let us take the $\NZ$-filtration $FR^+=\{F_qR^+\}_{q\in\NZ}$ of $R^+$ determined by $\B$, where
$$F_qR^+=K\hbox{-span}\{ u\in\B~|~d(u)\le q\},\quad q\in\NZ.$$
It is straightforward to check that $F_0R^+=K$, $F_qR^+\subseteq F_{q+1}R^+$ for all $q\in\NZ$, $R^+=\cup_{q\in\NZ}F_qR^+$, and, by referring to the proof of Theorem 2.3,
$$F_{q_1}R^+F_{q_2}R^+\subseteq F_{q_1+q_2}R^+,\quad q_1,q_2\in\NZ.$$
Hence, the filtration $FR^+$ constructed above makes $R^+$ into an $\NZ$-filtered algebra (indeed, one may check easily that the filtration $FR^+$ defined here coincides with the natural standard filtration of $R^+$, see also [Li6, Proposition A3.6]). Considering the associated $\NZ$-graded algebra $G(R^+)=\oplus_{q\in\NZ}G(R^+)_q$ of $R^+$, where $G(R^+)_0=K$, $G(R^+)_q=F_qR^+/F_{q-1}R^+$, $q\ge 1$, then $G(R^+)=K[\sigma(x_{ij})~|~(i,j)\in\Lambda_N]$, where $\sigma (x_{ij})$ is the homogeneous element in $G(R^+)$ represented by $x_{ij}$ and $d(\sigma (x_{ij}))=1=d(x_{ij})$, such that
$$\begin{array}{ll}\sigma (x_{mn})\sigma (x_{ij})=q^{-2}\sigma (x_{ij})\sigma (x_{mn}),&((i,j),(m,n))\in C_1\cup C_3,\\
\sigma (x_{mn})\sigma (x_{ij})=\sigma (x_{ij})\sigma (x_{mn}),&((i,j),(m,n))\in C_2\cup C_6,\\
\sigma (x_{mn})\sigma (x_{ij})=\sigma (x_{ij})\sigma (x_{mn})-(q^2-q^{-2})\sigma (x_{in})\sigma (x_{mj}),&((i,j),(m,n))\in
C_4,\\
\sigma (x_{mn})\sigma (x_{ij})=q^2\sigma (x_{ij})\sigma (x_{mn}),&((i,j),(m,n))\in C_5.\end{array}$$
Note that with the aid of the monomial ordering $\prec$ on $\B$ (as constructed in the proof of Theorem 2.3), we can further construct a graded monomial ordering $\prec_d$ on $\B$ as follows: for $u,v\in\B$,
$$u\prec_d\Leftrightarrow\left\{\begin{array}{l} d(u)<d(v),\\
\hbox{or} d(u)=d(v)~\hbox{and}~u\prec v.\end{array}\right.$$
It follows from [Li1, CH.IV, Theorem 4.1] that $G(R^+)$ is a quadratic solvable polynomial algebra with the PBW $K$-basis
$$\sigma(\B ) =\left\{\sigma (x_{i_1j_1})\sigma(x_{i_2j_2})\cdots
\sigma(x_{i_kj_k})~\left | ~\begin{array}{l}(i_{\ell},j_{\ell})\in\Lambda_N,~k\ge 1,\\
(i_1,j_1)\le_{lex} (i_2,j_2)\le_{lex}\cdots \le_{lex} (i_k,j_k)\end{array}\right.\right\}$$
and the graded monomial ordering $\prec_d$ induced by that on $R^+$. Noticing that $q^2\ne 0$, and $G(R^+)$ is an iterated skew polynomial algebra over the commutative polynomial ring $K[\sigma (x_{mn}),~\sigma (x_{ij})~|~((i,j),(m,n))\in C_2\cup C_6]$ (as one may see easily), it follows from[Lev],  [Li1, P. 176 Theorem 1.1], and the foregoing Proposition 3.2(ii) that $G(R^+)$ is an Auslander regular domain of global dimension $\frac{N(N+1)}{2}$ and $G(R^+)$ satisfies the Cohen-Macaulay property. Since the $\NZ$-filtration $FR^+$ of $R^+$ is a positive filtration and both $R^+$ and $G(R^+)$ are solvable polynomial algebras and hence are Noetherian domains, $R^+$ is then a Zariskian filtered ring in the sense of [LVO]. Also by [Li1] and [LVO] we know that  GK.dim$R^+=$ GK.dim$G(R^+)$, and GK.dim$M=$ GK.dim$G(M)$, $j_{R^+}(M)=j_{G(R^+)}(G(M))$ hold rue for every finitely generated $A$-module $M$, where $G(M)$ is the associated graded module of $M$ determined by a good filtration of $M$, and $j_{R^+}(M)$, respectively $j_{G(R^+)}(G(M))$, denotes the grade number of $M$, respectively the graded number of $G(M)$.   Therefore, it follows from [LVO, CH.III, Theorem 2, Theorem 6] that $R^+$ is an Auslander regular algebra satisfying the Cohen-Macaulay property.\par

(ii) That $K_0(R^+)\cong\mathbb{Z}$ follows from [Li1, P. 176, Theorem 1.1] (see also [LVO, P.125 Corollary 6.8], or [Li6, Theorem 3.3.3]).\par

This finishes the proof.\QED\v5

Also let us mention a result concerning modules over $R^+$ (resp. $R^-$).\v5

{\bf Theorem 3.5} Let the algebra $R^+$ (resp. $R^-$) be as before. Then the following statements hold.\par

(i) Let $L$ be a nonzero left ideal of $R^+$, and $R^+/L$  the left $R^+$-module. Considering Gelfand-Kirillov dimesion, we have GK.dim$R^+/L<$ GK.dim$A=\frac{N(N+1)}{2}$, and there is an algorithm for computing GK.dim$A/L$. Similar result holds true for $R^-$.\par

(ii)  Let $M$ be a finitely generated  $R^+$-module. Then a finite free resolution of $M$ can be algorithmically constructed, and the projective dimension of $M$ can be algorithmically computed. Similar result holds true for $R^-$.\vskip 6pt

{\bf Proof} (ii) That Gk.dim$R^+=\frac{N(N+1)}{2}$ follows from Theorem 3.3(ii). Since $R^+$ is a (quadric) solvable polynomial algebra by Theorem 2.3, it follows from [Li1, CH.V] that GK.dim$R^+/L<\frac{N(N+1)}{2}$ (this may also follow from classical Gelfand-Kirillov dimension theory [KL], for $R^+$ is now a Noetherian domain), and that there is an algorithm for computing GK.dim$R^+/L$.\par

(iii) This follows from [Li6, Ch.3].\QED\v5

We end this section by concluding that the algebra $R^+$ (resp. $R^-$) also has the elimination property for (one-sided) ideals in the sense of [Li5]  (see also [Li6, A3]), and that this property may be realized in a computational way. To see this clearly, let us first recall the Elimination Lemma given in [Li5].  Let $A=K[a_1,a_2\ldots ,a_n]$ be a  finitely generated $K$-algebra with
the PBW basis $\B =\{ a^{\alpha}=a_1^{\alpha_1}a_2^{\alpha_2}\cdots
a_n^{\alpha_n}~|~\alpha =(\alpha_1,\alpha_2\ldots
,\alpha_n)\in\NZ^n\}$ and, for a subset $U=\{ a_{i_1},a_{i_2},\ldots ,a_{i_r}\}\subset
\{ a_1,a_2,\ldots,a_n\}$ with $i_1<i_2<\cdots <i_r$,  let
$$T=\left\{ a_{i_1}^{\alpha_1}a_{i_2}^{\alpha_2}\cdots a_{i_r}^{\alpha_r}~\Big |~
(\alpha_1,\alpha_2,\ldots ,\alpha_r)\in\NZ^r\right\},\quad
V(T)=K\hbox{-span}T.$$\par

{\bf Lemma 3.6} ([Li5, Lemma 3.1]) Let the algebra $A$ and the notations be as fixed above, and let $L$ be a nonzero left ideal of $A$ and $A/L$ the left $A$-module defined by $L$.
If there is a subset  $U=\{ a_{i_1},a_{i_2}\ldots ,a_{i_r}\}
\subset\{a_1,a_2\ldots ,a_n\}$ with $i_1<i_2<\cdots <i_r$, such that $V(T)\cap L=\{ 0\}$, then $$\hbox{GK.dim}(A/L)\ge r.$$  Consequently, if  $A/L$ has finite GK dimension $\hbox{GK.dim}(A/L)=d<n$ ($=$
the number of generators of $A$), then  $$V(T)\cap L\ne
\{ 0\}$$ holds true for every subset $U=\{
a_{i_1},a_{i_2},\ldots ,a_{i_{d+1}}\}\subset$ $\{ a_1,a_2,\ldots ,a_n\}$ with
$i_1<i_2<\cdots <i_{d+1}$, in particular, for every $U=\{
a_1,a_2,\ldots ,a_s\}$ with $d+1\le s\le n-1$, we have $V(T)\cap
L\ne \{ 0\}$.\par\QED \v5

Once again, the discussion below will be carried out only for $R^+$, because similar argumentation works well for $R^-$.  Also for convenience of deriving the next theorem, let us write the set of generators of $R^+$ as $X=\{ x_1,x_2,\ldots x_{\omega}\}$ with $\omega =\frac{N(N+1)}{2}$, i.e., $R^+=K[x_1,x_2,\ldots ,x_{\omega}]$.
Thus, for a subset $U=\{ x_{i_1},x_{i_2},\ldots ,x_{i_r}\}\subset
\{ x_1, x_2,\ldots,x_{\omega}\}$ with $i_1<i_2<\cdots <i_r$, we write
$$T=\left\{ x_{i_1}^{\alpha_1}x_{i_2}^{\alpha_2}\cdots x_{i_r}^{\alpha_r}~\Big |~
(\alpha_1,\alpha_2,\ldots ,\alpha_r)\in\NZ^r\right\},\quad
V(T)=K\hbox{-span}T.$$\par

{\bf Theorem 3.7} With notation as fixed above, let $L$ be a left ideal of $R^+$. Then the following two statements hold.\par

(i) GK.dim$R^+/L< \omega$, and if GK.dim$R^+/L=d$, then
$$V(T)\cap L\ne\{ 0\}$$ holds true for every subset $U=\{
x_{i_1},x_{i_2},...,x_{d+1}\}\subset D$ with
$i_1<i_2<\cdots <i_{d+1}$, in particular, for every $U=\{
x_1,x_2\ldots x_s\}$ with $d+1\le s\le \omega -1$, we have $V(T)\cap
L\ne \{ 0\}$.\par

(ii) Without exactly knowing the numerical value GK.dim$R^+/L$, the elimination property for a left ideal $L=\sum_{i=1}^mR^+\xi_i$  of $R^+$ (resp. $R^-$) can be realized in a computational way, as follows:\par

Let $\prec$ be the monomial ordering on the PBW basis $\B$ of $R^+$ as constructed in the proof of Theorem 2.3 (or let $\prec_d$ be the graded monomial ordering as constructed in the proof of Theorem 3.4), and let $V(T)$ be as in (i). Then, employing  an elimination ordering  $\lessdot$ with respect to $T$ (which can always  be constructed if the existing monomial ordering on $\B$ is not an elimination ordering, see [Li6, Proposition 1.6.3]),  a Gr\"obner basis $\G$ of $L$ can be produced by running the noncommutative Buchberger algorithm for solvable polynomial algebras, such that
$$L\cap V(T)\ne \{ 0\} \Leftrightarrow \G\cap V(T)\ne \emptyset .$$\par

Similar statements as mentioned above hold true for $R^-$.\vskip 6pt

{\bf Proof} (i) Adopting the notations we fixed above, it follows from [Yam] that $R^+$ has the PBW basis $$\B =\{ x^{\alpha}=x_1^{\alpha_1}x_2^{\alpha_2}\cdots
x_{\omega}^{\alpha_{\omega}}~|~\alpha =(\alpha_1,\alpha_2\ldots
,\alpha_{\omega})\in\NZ^{\omega}\}.$$
Also by Theorem 3.3(ii) and Theorem 3.5(i) we have  GK.dim$R^+/L<\omega$. Therefore, the desired elimination property follows from  Lemma 3.6 mentioned above.\par

(ii) This follows from [Li6, Corollary 1.6.5].\QED\v5

{\bf Remark} Since $R^+$ (resp. $R^-$) is now a solvable polynomial algebra, if  $F=\oplus_{i=1}^sM_q(n)e_i$ is a free (left) $R^+$-module (resp. $R^-$-module) of finite rank, then a similar (even much stronger) result of Theorem 3.7 holds true for any finitely generated submodule $N=\sum_{i=1}^mM_q(n)\xi_i$ of $F$. The reader is referred to [Li6, Section 2.4] for a detailed argumentation.\v5

\centerline{Reference}{\parindent=.6truecm\par

\item{[BG]} ~R. Berger and V. Ginzburg, ~Higher ~symplectic
~reflection ~algebras and nonhomogeneous $N$-Koszul property, {\it
J. Alg.}, 1(304)(2006), 577--601.

\item{[Dr]} V.G. Drinfeld, Hopf algebras and the quantum Yang-Baxter equation, Doklady
Akademii Nauk SSSR 283(5) (1985) 1060--1064.

\item{[Jim]} M. Jimbo, A q-difference analogue of U(G) and the Yang-Baxter equation, Letters
in Mathematical Physics 10(1) (1985) 63--69.

\item{[KL]} G.R. Krause and T.H. Lenagan, {\it Growth of Algebras and
Gelfand-Kirillov Dimension}. Graduate Studies in Mathematics.
American Mathematical Society, 1991.

\item{[K-RW]} A. Kandri-Rody and V. Weispfenning, Non-commutative
Gr\"obner bases in algebras of solvable type. {\it J. Symbolic
Comput.}, 9(1990), 1--26. Also available as: Technical Report University of Passau,
MIP-8807, March 1988.

\item{[Lev]} T. Levasseur, Some properties of noncommutative regular graded rings.
{\it Glasgow Math. J}., 34(1992), 277--300.

\item{[Li1]} H. Li, {\it Noncommutative Gr\"obner Bases and Filtered-graded Transfer}.
Lecture Notes in Mathematics, Vol. 1795, Springer, 2002.

\item{[Li2]} H. Li, $\Gamma$-leading ~homogeneous~ algebras~ and
Gr\"obner bases. In: {\it Recent Developments in Algebra and
Related Areas} (F. Li and C. Dong eds.), Advanced Lectures in
Mathematics, Vol. 8, International Press \& Higher Education Press,
Boston-Beijing, 2009, 155 -- 200. \url{arXiv:math.RA/0609583,
http://arXiv.org}

\item{[Li3]} H. Li, {\it Gr\"obner Bases in Ring Theory}. World Scientific Publishing Co., 2011. \url{https://doi.org/10.1142/8223}\par

\item{[Li4]} H. Li, A note on solvable polynomial algebras. {\it Computer Science Journal of Moldova}, vol.22, no.1(64), 2014, 99--109.  arXiv:1212.5988 [math.RA]

\item{[Li5]} H. Li, An elimination lemma for algebras with PBW bases. {\it Communications in
Algebra}, 46(8)(2018), 3520-3532.\par

\item{[Li6]} H. Li, {\it Noncommutative polynomial algebras of solvable type and their modules: Basic constructive-computational theory and methods}. Chapman and Hall/CRC Press, 2021.

\item{[LVO]} H. Li, F. Van Oystaeyen,  {\it Zariskian Filtrations}. K-Monograph
in Mathematics, Vol.2. Kluwer Academic Publishers, 1996; Berlin Heidelberg: Springer-Verlag, 2003.

\item{[LS]} V. Levandovskyy and H. Sch\"onemann, Plural: a computer algebra system for noncommutative polynomial algebras. In: {\it Proc. Symbolic and Algebraic Computation}, International Symposium ISSAC 2003, Philadelphia, USA, 176--183, 2003.

\item{[Mor]} T. Mora, An introduction to commutative and
noncommutative Gr\"obner Bases. {\it Theoretic Computer Science},
134(1994), 131--173.

\item{[Pos]} ~L. Positselski, Nonhomogeneous quadratic duality and
curvature, {\it Funct. Anal. Appl.}, 3(27)(1993), 197--204.

\item{[Ros]} M. Rosso, Finite dimensional representations of the quantum analogue of the
enveloping algebra of a complex simple Lie algebra, Comm. Math. Phys. 117
(1988) 581--593.

\item{[Yam]} I. Yamane, A Poincare-Birkhoff-Witt theorem for
quantized universal enveloping algebras of type $A_N$, Publ., {\it
RIMS. Kyoto Univ}., 25(3)(1989), 503--520.}

\end{document}